\documentclass[a4paper,12pt]{article}
\usepackage{comment}
\usepackage{cite}
\usepackage{amsmath}
\usepackage{amssymb}
\usepackage{amsfonts}
\usepackage[T1]{fontenc}
\usepackage[utf8]{inputenc}
\usepackage{graphicx}
\usepackage{fancyhdr}
\usepackage{float}
\usepackage{xcolor}
\usepackage{authblk}
\usepackage{empheq}
\usepackage[hyphens]{url}
\usepackage{hyperref} 
\usepackage[]{breakurl}
\usepackage{mathrsfs}
\usepackage[]{amsthm}

\usepackage{geometry}
\begin{document}

\title{Asymptotic equivalents of partial sums of the reciprocals of prime numbers via the von Mangoldt function}

\author[$\dagger$]{Jean-Christophe {\sc Pain}$^{1,2,}$\footnote{jean-christophe.pain@cea.fr}\\
\small
$^1$CEA, DAM, DIF, F-91297 Arpajon, France\\
$^2$Universit\'e Paris-Saclay, CEA, Laboratoire Mati\`ere en Conditions Extr\^emes,\\ 
F-91680 Bruy\`eres-le-Ch\^atel, France
}

\date{}

\maketitle

\begin{abstract}
In this paper, we discuss an alternative approach to determine an asymptotic equivalent of the partial sum of the reciprocals of prime numbers. This well-known result, related to Mertens' second theorem, is usually derived through methods similar to those found in Hardy and Wright's book ``An introduction to the theory of numbers'', involving comparisons with integrals. The present proof differs in several respects, combining an equivalent for the partial sum of $\Lambda(m)/m$, where $\Lambda$ denotes the von Mangoldt function, with an application of Abel's summation formula and properties of the second Chebyshev function $\Psi(x)=\sum_{n\le x}\Lambda(n)$. A simple application to the study of integers with large prime factors is also presented. Beyond the pedagogical aspect of this work, the aim is to highlight the complementarity of arithmetic functions and to show that interesting (and nontrivial) results can be obtained by means of elementary methods.

\end{abstract}

\section{Introduction}\label{sec1}

The divergence of the sum of the reciprocals of all prime numbers
\begin{equation*}
    \sum _{p{\text{ prime}}}{\frac {1}{p}}=\infty,
\end{equation*}
was proved by Euler in 1737. There are many different proofs of this divergence. For instance, from Dusart's inequality \cite{Dusart1999}:
\begin{equation*}
    p_{n}<n\log n+n\log \log n
\end{equation*}
for $n\geq 6$, we get immediately
\begin{align*}
    \sum _{n=1}^{\infty }{\frac {1}{p_{n}}}&\geq \sum _{n=6}^{\infty }{\frac {1}{p_{n}}}\\&\geq \sum _{n=6}^{\infty }{\frac {1}{n\log n+n\log \log n}}\\&\geq \sum _{n=6}^{\infty }{\frac {1}{2n\log n}}=\infty. 
\end{align*}
It is also possible to use a lower bound for the partial sums stating that 
\begin{equation*}
    \sum _{\scriptstyle p{\text{ prime}} \atop \scriptstyle p\leq n}{\frac {1}{p}}\geq \log \log(n+1)-\log\left(\frac {\pi ^{2}}{6}\right).
\end{equation*}
Due to the occurrence of the double natural logarithm, the divergence is slow. The study of the partial sums of the reciprocals of prime numbers, aa well as of the reciprocals of the largest prime factor of an integer, plays an important role in number theory \cite{Alladi1979,Ivic1981,Erdos1986,Koninck1993,Koninck1994}.  

The usual way of determining an asymptotic equivalent of the partial sum of the reciprocals of prime numbers proceeds as follows. We shall show that, as \(x\to\infty\),
\begin{equation*}
    \sum_{p\le x}\frac{1}{p}=\log\log x + b + O\!\Big(\frac{1}{\log x}\Big),
\end{equation*}
where the constant \(b\) is defined by the limit
\begin{equation*}
    b=\lim_{x\to\infty}\Big(\sum_{p\le x}\frac{1}{p}-\log\log x\Big).
\end{equation*}
Let us define
\begin{equation*}
    S(x):=\sum_{p\le x}\frac{1}{p},
\end{equation*}
and express \(S(x)\) as a Stieltjes integral with respect to the function \(\pi(t)\) (the number of primes \(\le t\)):
\begin{equation*}
    S(x)=\int_{2^-}^{x}\frac{1}{t}\,d\pi(t).
\end{equation*}
Applying integration by parts for Stieltjes integrals (or the discrete version, summation by parts), we obtain
\begin{equation*}
    \int_{2^-}^{x}\frac{1}{t}\,d\pi(t) = \frac{\pi(x)}{x} + \int_{2}^{x}\frac{\pi(t)}{t^{2}}\,dt.
\end{equation*}
We now use the Prime Number Theorem in the form
\begin{equation*}
    \pi(t)=\frac{t}{\log t}+O\!\Big(\frac{t}{\log^{2}t}\Big)\qquad(t\ge 3).
\end{equation*}
Substituting this into the previous expression, we get for the first term:
\begin{equation*}
    \frac{\pi(x)}{x}=\frac{1}{\log x}+O\!\Big(\frac{1}{\log^{2}x}\Big).
\end{equation*}
For the integral,
\begin{equation*}
\int_{2}^{x}\frac{\pi(t)}{t^{2}}\,dt
=\int_{2}^{x}\frac{1}{t\log t}\,dt
+\int_{2}^{x}O\!\Big(\frac{1}{t\log^{2}t}\Big)\,dt.
\end{equation*}
We now study these two integrals separately. The first integral is
\begin{equation*}
  \int_{2}^{x}\frac{1}{t\log t}\,dt = \log\log x - \log\log 2 .
\end{equation*}
For the second, uniformly for \(x\ge 3\),
\begin{equation*}
  \int_{2}^{x}O\!\Big(\frac{1}{t\log^{2}t}\Big)\,dt
  = O\!\Big(\int_{2}^{x}\frac{dt}{t\log^{2}t}\Big)
  = O\!\Big(\frac{1}{\log x}\Big),
\end{equation*}
since 
\begin{equation*}
\int_{2}^{x}\frac{dt}{t\log^{2}t}=\Big[-\frac{1}{\log t}\Big]_{2}^{x}= \frac{1}{\log 2}-\frac{1}{\log x}=O\!\Big(\frac{1}{\log x}\Big).
\end{equation*}
Gathering all contributions, we get
\begin{equation*}
    S(x) = \frac{1}{\log x}+ \big(\log\log x-\log\log 2\big) + O\!\Big(\frac{1}{\log x}\Big).
\end{equation*}
The term \(\frac{1}{\log x}\), combined with the error term \(O(1/\log x)\), still gives an \(O(1/\log x)\) contribution. Hence, there exists a constant \(b\) (equal to \(-\log\log 2\) plus constants arising from the error terms and the contribution of small primes) such that
\begin{equation*}
    \sum_{p\le x}\frac{1}{p}=\log\log x + b + O\!\Big(\frac{1}{\log x}\Big),
\end{equation*}
which completes the proof. The constant \(b\) can be expressed explicitly in terms of Euler’s constant \(\gamma\) and a series over the primes:
\begin{equation*}
    b=\gamma + \sum_{p}\Big(\log\big(1-\tfrac{1}{p}\big)+\tfrac{1}{p}\Big),
\end{equation*}
which follows from the analysis of the Euler product for \(\zeta(s)\) and the behavior of \(\log\zeta(s)\) near \(s=1\).  
The proof above avoids these calculations and merely establishes the existence of the constant \(b\) and the \(O(1/\log x)\) error term.

The main purpose of the present work is to highlight the differences and connections between the above derivation, and an other one, based on elementary manipulations of arithmetic functions, and explained below.

For $n \ge 1$, the von Mangoldt function reads

\begin{equation}\label{vonman}
    \Lambda(n)=
    \begin{cases}
    \log p &\text{if }n=p^{\alpha}\ (\alpha\ge1),\\[4pt]
    0 &\text{otherwise.}
    \end{cases}
\end{equation}
We will use the asymptotic equivalent \cite{Apostol1976} (see section \ref{sec2}): 
\begin{equation*}
    \sum_{n\le x}\frac{\Lambda(n)}{n} = \log x + O(1),
\end{equation*}
to recover (see section \ref{sec3}):
\begin{equation*}
    \sum_{p\le x}\frac{1}{p} = \log\log x + b + O\!\Big(\frac{1}{\log x}\Big).
\end{equation*}
and the derivation will be applied to the determination of the density of integers with large prime factors in section \ref{sec4}.

\section{Asymptotic formula involving the von Mangoldt function}\label{sec2}

\subsection{Expression of $\log(n!)$ in terms of $\Lambda$ using the Legendre formula}

The Legendre formula gives, for every prime number $p$:
\begin{equation*}
    \nu_p(n!) = \sum_{k\ge1}\Big\lfloor\frac{n}{p^k}\Big\rfloor,
\end{equation*}
where $\lfloor x \rfloor$ denotes the integer part of $x$ and $\nu_p(\ell)$ the $p-$adic valuation of $\ell$. We thus have, writing \(\log(n!)\) as a sum over the contributions of the different prime factors:
\begin{equation*}   
    \log(n!)=\sum_{p\ \text{prime}} v_p(n!)\log p = \sum_{p}\sum_{k\ge1}\Big\lfloor\frac{n}{p^k}\Big\rfloor\log p.
\end{equation*}
According to the definition of the von Mangoldt function (\ref{vonman}), the double sum above can be written as the sum over all prime powers \(p^k\le n\) :
\begin{equation*}  
    \log(n!)=\sum_{p}\sum_{k\ge1}\Big\lfloor\frac{n}{p^k}\Big\rfloor\Lambda(p^k)=\sum_{m\le n}\Big\lfloor\frac{n}{m}\Big\rfloor\Lambda(m),
\end{equation*}
since only the \(m\) of the kind \(p^k\) contribute (the others have \(\Lambda(m)=0\)). Thus we recover
\begin{equation}   
    \log(n!)=\sum_{m\le n}\Lambda(m)\Big\lfloor\frac{n}{m}\Big\rfloor.
\end{equation}

\subsection{A more direct proof of $\log(n!)=\sum_{m\le n}\Lambda(m)\Big\lfloor\frac{n}{m}\Big\rfloor$}

Let us start from
\begin{equation}
    \log(n!) = \sum_{k=1}^{n} \log k,
\end{equation}
and then use the fundamental property of the von Mangoldt function:
\begin{equation*}  
    \log k = \sum_{d \mid k} \Lambda(d), \qquad \text{for all } k \ge 1,
\end{equation*}
since if \( k = \prod_{p} p^{a_p} \), then the only divisors \( d \) such that \( \Lambda(d) \neq 0 \) are the prime powers \( p^j \), and
\begin{equation*}  
    \sum_{j=1}^{a_p} \Lambda(p^j) = a_p \log p,
\end{equation*}
which gives
\begin{equation*}  
    \sum_{d \mid k} \Lambda(d) = \sum_p a_p \log p = \log k.
\end{equation*}
Substituting the latter identity in the sum (\ref{inter}), we get
\begin{equation*}  
    \log(n!) = \sum_{k=1}^{n} \sum_{d \mid k} \Lambda(d).
\end{equation*}

Let us now change the order of the summations. For each \( d \le n \), the number of integers \( k \le n \) such that \( d \mid k \) is given by \( \left\lfloor \tfrac{n}{d} \right\rfloor \). Thus
\begin{equation*}  
\log(n!) = \sum_{d \le n} \Lambda(d) \, \Big\lfloor \frac{n}{d} \Big\rfloor.
\end{equation*}
Renaming the variable \( d \) into \( m \), we obtain:
\begin{equation}\label{inter}  
\log(n!) = \sum_{m \le n} \Lambda(m)\, \Big\lfloor \frac{n}{m} \Big\rfloor.
\end{equation}

\subsection{On the sum $\sum_{m\le x}\frac{\Lambda(m)}{m}$}

Let $x$ be a natural number. Writing $\lfloor x/m\rfloor = x/m + O(1)$, we get, from Eq. (\ref{inter}):
\begin{equation*}
    \log(x!) = x\sum_{m\le x}\frac{\Lambda(m)}{m} + O\!\Big(\sum_{m\le x}\Lambda(m)\Big) = x\sum_{m\le x}\frac{\Lambda(m)}{m} + O(\Psi(x)),
\end{equation*}
where
\begin{equation*}
    \Psi(x)=\sum_{n\le x}\Lambda(n)
\end{equation*}
represents the second Chebyshev function. We know (see Appendix A) that $\Psi(x)=O(x)$ thanks to an elementary argument. Therefore,
\begin{equation*}
    \frac{\log(x!)}{x} = \sum_{m\le x}\frac{\Lambda(m)}{m} + O(1).
\end{equation*}
By the Stirling formula,
\begin{equation*}
    \log(x!) = x\log x - x + O(\log x),
\end{equation*}
we get
\begin{equation*}
    \frac{\log(x!)}{x} = \log x - 1 + O\!\Big(\frac{\log x}{x}\Big) = \log x + O(1),
\end{equation*}
and finally,
\begin{equation}\label{newman}
    \sum_{m\le x}\frac{\Lambda(m)}{m} = \log x + O(1).
\end{equation}
According to M\"obius inversion, the latter formula is equivalent to
\begin{equation*}
    \sum_{m\le x}\frac{1}{m}\sum_{d\mid m}\mu(d)\log\left(\frac{m}{d}\right) = \log x + O(1).
\end{equation*}

\section{Equivalent of $\sum_{p\le x}\frac{1}{p}$}\label{sec3}

We first gather the terms of the previous sum according to powers of primes:
\begin{align*}
    \sum_{m\le x}\frac{\Lambda(m)}{m}
    &= \sum_{p\le x}\sum_{\substack{k\ge1\\p^k\le x}}\frac{\log p}{p^k}\nonumber\\
    &= \sum_{p\le x}\frac{\log p}{p}\Big(1+\frac{1}{p}+\frac{1}{p^2}+\cdots\Big) + O\!\Big(\sum_{p\le x}\frac{\log p}{p^2}\Big).
\end{align*}

The series 
\begin{equation*}
    \sum_p \frac{\log p}{p^2}
\end{equation*}
converges, so the error term is $O(1)$. Thus, from Eq. (\ref{newman}), we get
\begin{equation}\label{mertens1}
    \sum_{p\le x}\frac{\log p}{p} = \log x + O(1),
\end{equation}
which is related to Mertens' first theorem \cite{Mertens1874,Mertens1874b,Tenenbaum1995,OEIS}, stating that
\begin{equation*}
    \sum_{p\leq n}\frac{\log p}{p}-\log n
\end{equation*}
does not exceed 2 in absolute value for any $n\geq 2$.

Let us now set
\begin{equation*}
    \mathscr{A}(t)=\sum_{p\le t}\frac{\log p}{p}
\end{equation*}
and apply Abel's summation formula to the sequence
\begin{equation*}
    a_p=\frac{\log p}{p}
\end{equation*}
and the function 
\begin{equation*}
    f(t)=\frac{1}{\log t}.
\end{equation*}
This yields
\begin{equation*}
    \sum_{p\le x}\frac{1}{p} = \mathscr{A}(x)f(x) - \int_2^x \mathscr{A}(t) f'(t)\,dt.
\end{equation*}
Since $f'(t)=-1/(t(\log t)^2)$ and, from (\ref{mertens1}), $\mathscr{A}(t)=\log t+O(1)$, we get
\begin{align*}
    \sum_{p\le x}\frac{1}{p} &= \frac{\mathscr{A}(x)}{\log x} + \int_2^x \frac{\mathscr{A}(t)}{t(\log t)^2}\,dt\nonumber\\
                             &= \frac{\log x + O(1)}{\log x} + \int_2^x\frac{\log t + O(1)}{t(\log t)^2}\,dt.
\end{align*}
The first term equals $1+O(1/\log x)$, and
\begin{align*}
    \int_2^x\frac{\log t}{t(\log t)^2}\,dt &= \int_2^x\frac{1}{t\log t}\,dt\nonumber\\
                                           &= \log\log x + C_1.
\end{align*}
The integral of the error term gives another constant $O(1)$. Hence
\begin{equation}\label{sumpar}
    \sum_{p\le x}\frac{1}{p} = \log\log x + b + O\!\Big(\frac{1}{\log x}\Big).
\end{equation}

\section{Application to the study of integers with large prime factors}\label{sec4}

For each natural number $n\geq 2$, let us consider that $n$ has large prime factors, if the largest prime factor entering its decomposition is larger than $\sqrt{n}$. We denote by $\mathscr{N}(x)$ the ensemble of natural numbers satisfying that property. According to Eq. (\ref{sumpar}), we have
\begin{align}\label{prem}
    \sum_{\sqrt{n}<p\le n}\frac{1}{p}&=\sum_{p\le n}\frac{1}{p}-\sum_{p\le\sqrt{n}}\frac{1}{p}\nonumber\\
                                     &= \log(\log n)+b+o(1)-\bigl(\log(\log\sqrt{n})+b+o(1)\bigr)\nonumber\\
                                     &= \log(\log n)-\log\!\Bigl(\frac{\log n}{2}\Bigr)+o(1)\nonumber\\
                                     &= \log(\log n)-\log(\log n)+\log 2+o(1),
\end{align}
which implies that the sequence indeed converges and its limit is $\log 2$.

Let us first assume that $n = qp$ is a number $\le x$ having a prime factor strictly larger than $\sqrt{n}$. Since $p$ is the largest prime factor of $n$, we must have $p > \sqrt{n}$. Then $q = n/p$ is strictly smaller than $n/\sqrt{n} = \sqrt{n} < p$. Moreover, since $n \le x$, we have 
\begin{equation*}
    p = n/q \le x/q. 
\end{equation*}
Conversely, suppose that $q < p \le x/q$. Then $\sqrt{n} = \sqrt{qp} < \sqrt{p}\sqrt{p} = p$, so $n$ indeed has a large prime factor. If $q$ is another prime factor of $n$, it divides $n/p$, and $n/p < n/\sqrt{n} = \sqrt{n}$; thus $p$ is indeed the largest prime factor of $n$. Finally, 
\begin{equation*}
    n = qp \le qx/q = x.
\end{equation*}
Let us assume that $qp = q'p'$. By the previous question, both $p$ and $p'$ are the largest prime factors of $qp = q'p'$. Hence $p = p'$, and it follows that 
\begin{equation*}
    q = (qp)/p = (q'p')/p' = q'. 
\end{equation*}
The converse is obvious.

Subsequently, the integers belonging to $\mathscr{N}(x)$ are precisely those of the form $qp$ where $p$ is prime and $q$ is a positive integer such that $q < p \le x/q$. Moreover, each such integer corresponds to a unique pair $(p,q)$ with $p$ prime and $q$ a positive integer satisfying 
\begin{equation}\label{ineg}
    q < p \le x/q.
\end{equation}
The number $\mathscr{G}(x)$ of pairs $(p,q)$ with $p$ prime and $q$ an integer satisfying (\ref{ineg}). For a fixed prime $p$, an integer $q$ satisfies the condition if and only if $q < p$ and $q \le x/p$, that is, if and only if
\begin{equation*}
    q \le \min\!\Bigl(p-1,\ \Bigl\lfloor \frac{x}{p}\Bigr\rfloor\Bigr).
\end{equation*}
For a given prime $p$, the number of integers $n\in \mathscr{N}(x)$ that can be written $n=qp$ with $q < p \le x/q$ is therefore
\begin{equation*}
    \min\!\Bigl(p-1,\ \Bigl\lfloor \frac{x}{p}\Bigr\rfloor\Bigr).
\end{equation*}
Possible primes $p$ being those $\le x$, one has
\begin{equation*}
    \mathscr{G}(x)=\sum_{p\le x}\min\!\Bigl(p-1,\ \Bigl\lfloor \frac{x}{p}\Bigr\rfloor\Bigr).
\end{equation*}
Since $p-1$ is an integer, we have the equivalence
\begin{equation*}
    p-1 \le \Bigl\lfloor\frac{x}{p}\Bigr\rfloor
\iff p-1 \le \frac{x}{p}.
\end{equation*}
The right-hand inequality is equivalent (for $p\ge0$) to $p^2 - p \le x$,
that is,
\begin{equation*}
    (p-\tfrac12)^2 \le x+\tfrac14,
\end{equation*}
and since $p-\tfrac12>0$, this is equivalent to
\begin{equation*}
    p-\tfrac12 \le \sqrt{x+\tfrac14}\iff p \le \tfrac12+\sqrt{\tfrac14+x} =: \eta(x).
\end{equation*}
Since $x>0$, we have
\begin{equation*}
    \sqrt{x}<\sqrt{\tfrac14+x}<\tfrac12+\sqrt{\tfrac14+x}=\eta(x).
\end{equation*}
As $x\ge1$, we have $x\le x^2$,
\begin{equation*}
    \eta(x)=\tfrac12+\sqrt{\tfrac14+x}\le\tfrac12+\sqrt{\tfrac14+x^2}.
\end{equation*}
By the triangle inequality, one gets 
\begin{equation*}
    \eta(x)\le\tfrac12+\sqrt{\tfrac14}+\sqrt{x^2}=\tfrac12+\tfrac12+x=1+x.
\end{equation*}
From the above results, we have
\begin{align*}
    \mathscr{G}(x)&=\sum_{p\le x}\min\!\Bigl(p-1,\Bigl\lfloor\frac{x}{p}\Bigr\rfloor\Bigr)\nonumber\\
        &=\sum_{p\le\eta(x)}(p-1)+\sum_{p>\eta(x)}\Bigl\lfloor\frac{x}{p}\Bigr\rfloor.
\end{align*}
We can rewrite this as
\begin{equation*}
    \mathscr{G}(x)=\sum_{p\le\sqrt{x}}(p-1)+\sum_{\sqrt{x}<p\le\eta(x)}(p-1)+\sum_{p>\eta(x)}\Bigl\lfloor\frac{x}{p}\Bigr\rfloor.
\end{equation*}
Suppose there is no prime in the interval $(\sqrt{x},\eta(x)]$. Then the second sum contributes zero, and the third sum equals
\begin{equation*}
    \sum_{\sqrt{x}<p\le x}\Bigl\lfloor\frac{x}{p}\Bigr\rfloor,
\end{equation*}
as required.

Now suppose there exists a prime $p$ in $(\sqrt{x},\eta(x)]$. This interval having length $<1$, there is at most one such prime. Then $p\le\eta(x)$, implying 
\begin{equation*}
    p-1\le\lfloor x/p\rfloor.
\end{equation*}
In addition, $\sqrt{x}<p$ implies $x<p^2$ and hence $x/p<p$, so $\lfloor x/p\rfloor<p$. From $p-1\le\lfloor x/p\rfloor<p$ we deduce $\lfloor x/p\rfloor=p-1$. Thus, in this case, the sum of the second and third sums equals
\begin{equation*}
    \sum_{\sqrt{x}<p\le x}\Bigl\lfloor\frac{x}{p}\Bigr\rfloor.
\end{equation*}
We have
\begin{equation*}
    \sum_{p\le\sqrt{x}}(p-1)\le \sum_{p\le\sqrt{x}}\sqrt{x} = \pi(\sqrt{x})\sqrt{x}.
\end{equation*}
If $\pi(x)$ is the usual prime counting function defined in the introduction, we can bound $\pi(n)\le e\,n/\log n$ (see Appendix B), yielding
\begin{equation*}
    \sum_{p\le\sqrt{x}}(p-1)\le \pi(\sqrt{x})\sqrt{x}\le e\frac{\sqrt{x}}{\log\sqrt{x}}= e\,\frac{x}{\log\sqrt{x}}= o(x).
\end{equation*}
From Eq. (\ref{prem}), we have
\begin{equation*}
    \sum_{\sqrt{x}<p\le x}\Bigl(\frac{x}{p}-\Bigl\lfloor\frac{x}{p}\Bigr\rfloor\Bigr)\le \sum_{\sqrt{x}<p\le x}\frac{1}{p}= \log 2 + o(1).
\end{equation*}
Hence
\begin{align*}
    \sum_{\sqrt{x}<p\le x}\Bigl\lfloor\frac{x}{p}\Bigr\rfloor&=\sum_{\sqrt{x}<p\le x}\frac{x}{p}+o(x)\nonumber\\
                                                             &= x\log 2 + o(x),
\end{align*}
where we used again Eq. (\ref{prem}). Combining all the above results, we get
\begin{align*}
\mathscr{G}(x)&=\sum_{p\le\sqrt{x}}(p-1)+\sum_{\sqrt{x}<p\le x}\Bigl\lfloor\frac{x}{p}\Bigr\rfloor\nonumber\\
              &= o(x)+x\log 2+o(x)\nonumber\\
              &= x\log 2+o(x).
\end{align*}
It follows that
\begin{align*}
\frac{\mathscr{G}(n)}{n}&=\frac{n\log 2+o(n)}{n}\nonumber\\
              &=\log 2+o(1)
\end{align*}
has a limit equal to $\log 2$, which is the density of the set of integers having large prime factors has  $\log 2$.

\section{Conclusion}

We have detailed an alternative proof of the asymptotic equivalent of partial sums of reciprocals of prime numbers. The proof is elementary, in the sense that it does not rely on the prime number theorem, requiring the asymptotic behavior of the prime counting function $\pi(x)$. The main idea consists in resorting to the von Mangoldt arithmetic function and to proove the identity
\begin{equation*}
    \sum_{n\le x}\frac{\Lambda(n)}{n} = \log x + O(1).
\end{equation*}
The next steps of the derivation (more classical) consists in using the corollary of the first Mertens' theorem
\begin{equation*}
    \sum_{p\le x}\frac{\log p}{p} = \log x + O(1),
\end{equation*}
and, performing Abel transforms, to obtain, with the Chebyshev bound $\Psi(x)=O(x)$, the usual equivalent
\begin{equation*}
    \sum_{p\le x}\frac{1}{p} = \log\log x + b + O\!\Big(\frac{1}{\log x}\Big).
\end{equation*}
This last part of the derivation is common to the proof of the prime number theorem, which we formulated in terms of the Stieltjes integral. 

Let us recall that such a result is strongly related to the second Mertens theorem, which states that, if $x$ is a positive real number, one has \cite{Hardy1954}:
\begin{equation*}
    M=\lim _{n\rightarrow \infty }\left(\sum _{\scriptstyle p{\text{ prime}} \atop \scriptstyle p\leq n}{\frac {1}{p}}-\log(\log n)\right)=\gamma +\sum _{p}\left[\log \!\left(1-{\frac {1}{p}}\right)+{\frac {1}{p}}\right],
\end{equation*}
where $M$ is the Meissel-Mertens constant (see for instance Ref. \cite{Finch2003,Havil2019}):
\begin{equation*}
    M=\gamma+\sum_{p}^{+\infty}\left(\log \left(1-\frac {1}{p}\right)+\frac {1}{p}\right),
\end{equation*}
and $\gamma$ the Euler-Mascheroni constant. One has $M\approx 0.2614972128$.

In addition, the relation ($\zeta(s)$ is the Riemann zeta function and $\text{Re}(s)>1$):
\begin{equation*}
    \log \zeta (s)=\sum _{n=2}^{\infty }\frac{\Lambda (n)}{\log(n)}\,\frac {1}{n^{s}},
\end{equation*}
which plays a major role in the theory of Dirichlet series, may be useful to derive asymptotic equivalents of partial sums.

Moreover, in the same framework, it would be worth taking advantage of the Selberg identity \cite{Selberg1949}:
\begin{equation*}
\Lambda (n)\log(n)+\sum _{d\,|\,n}\Lambda (d)\Lambda \!\left({\frac {n}{d}}\right)=\sum _{d\,|\,n}\mu (d)\log ^{2}\left({\frac {n}{d}}\right)
\end{equation*}
as well as of generalized von Mangoldt function 
\begin{equation*}
    \Lambda _{k}(n)=\sum \limits _{d\mid n}\mu (d)\log ^{k}\left(\frac{n}{d}\right)
\end{equation*}
where $k$ is a positive integer.

\section*{Appendix A: Bounds and equivalent of the second Chebyshev function}\label{sec:chebyshev}

In this appendix we recall the well-known bounds of the Chebyshev function $\psi$ in order to show that $\Psi(x)=O(x)$. For that purpose, we first use the binomial coefficient
\begin{equation*}
    \binom{2n}{n} = \frac{(2n)!}{(n!)^2}.
\end{equation*}
On one hand, according the the Newton binomial expansion: 
\begin{equation*}
    (1+1)^{2n}=\sum_{k=0}^{2n}\binom{2n}{k}\geq\binom{2n}{n}
\end{equation*}
and thus
\begin{equation}\label{combin}
    \binom{2n}{n} \le 4^n. 
\end{equation}
On the other hand, we can express this quantity in terms of prime numbers. For each prime $p$,
\begin{equation*}
    \nu_p(m!) = \sum_{k\ge1}\Big\lfloor\frac{m}{p^k}\Big\rfloor\quad\text{and thus}\quad\nu_p\left(\binom{2n}{n}\right)
              = \sum_{k\ge1}\!\left(\Big\lfloor\frac{2n}{p^k}\Big\rfloor - 2\Big\lfloor\frac{n}{p^k}\Big\rfloor\right).
\end{equation*}
Each term is $0$ or $1$, hence $v_p(\binom{2n}{n})\ge1$ whenever $p\in(n,2n]$. Taking logarithms:
\begin{equation*}
    \log\binom{2n}{n} = \sum_{p\le 2n} \nu_p\left(\binom{2n}{n}\right)\log p.
\end{equation*}
Therefore,
\begin{equation*}
    \sum_{n<p\le 2n}\log p \;\le\; \log\binom{2n}{n} \;\le\; \log((2n)!).
\end{equation*}
The combinatorial inequality (\ref{combin}) gives
\begin{equation*}
    \log\binom{2n}{n} \le 2n\log 2,
\end{equation*}
so that
\begin{equation*}
    \sum_{n<p\le 2n}\log p \le 2n\log 2.
\end{equation*}
Letting $\Psi(x)=\sum_{k\in\mathbb{N}}\sum_{p^k\le x}\log p$, we obtain
\begin{equation*}
    \Psi(2n)-\Psi(n)\le 2n\log 2,
\end{equation*}
and summing this inequality over dyadic intervals $[2^k,2^{k+1}]$ gives $\Psi(x)=O(x)$.

We also have $\displaystyle \binom{2n}{n}\ge \frac{4^n}{2n+1}$, hence
\begin{equation*}
    \log\binom{2n}{n}\ge 2n\log 2 - O(\log n),
\end{equation*}
and therefore
\begin{equation*}
    \sum_{n<p\le 2n}\log p \ge 2n\log 2 - O(\log n).
\end{equation*}
We deduce the existence of constants $c_1,c_2>0$ such that
\begin{equation*}
    c_1x \le \Psi(x)\le c_2x \quad (x\ge1),
\end{equation*}
that is,
\begin{equation*}
    \Psi(x)=O(x).
\end{equation*}

\section*{Appendix B: Upper bound for the prime counting function}

Let us consider two natural numbers $n_1$ and $n_2$ such that $0<n_2/2\leq n_1<n_2$. Let $p$ be a prime number such that $n_1 < p \le n_2$. Since $n_1 < p$, we have $v_p(n_1!) = 0$, and since $n_2 - n_1 \le n_2 - n_2/2 = n_2/2 \le n_1 < n_2$, we also have $v_p((n_2-n_1)!) = 0$. It follows that
\begin{equation*}
    \nu_p\!\left(\binom{n_2}{n_1}\right) = v_p(n_2!) \ge 1.
\end{equation*}
Hence, the prime number $p$ divides the binomial coefficient $\binom{n_2}{n_1}$.  
As $p$ is any prime satisfying $a < p \le b$, the product
\begin{equation*}
    \prod_{n_1 < p \le n_2} p
\end{equation*}
divides $\binom{n_2}{n_1}$. For $n_1 = m + 1$ and $n_2 = 2m + 1$, we have
\begin{equation*}
    \frac{n_2}{2} = m + \frac{1}{2} \le n_1 = m + 1 < 2m + 1 = n_2.
\end{equation*}
Hence, the integer
\begin{equation*}
    \prod_{m + 1 < p \le 2m + 1} p
\end{equation*}
divides
\begin{equation*}
    \binom{2m + 1}{m + 1}=\binom{2m + 1}{m}.
\end{equation*}
We have
\begin{equation*}
    2^{2m + 1} = (1 + 1)^{2m + 1} = \sum_{k = 0}^{2m + 1} \binom{2m + 1}{k} \ge \binom{2m + 1}{m} + \binom{2m + 1}{m + 1} = 2\binom{2m + 1}{m}.
\end{equation*}
Dividing both sides by $2$ gives
\begin{equation*}
    2^{2m} \ge \binom{2m + 1}{m}.
\end{equation*}
We therefore obtain
\begin{equation*}   
    \prod_{m + 1 < p \le 2m + 1} p\le \binom{2m + 1}{m + 1} = \binom{2m + 1}{m} \le 4^m.
\end{equation*}
For $n \ge 1$, let $\mathscr{P}_n$ denote the property:
\begin{equation*}
    \text{For all } k \in \{1, \dots, 2n\},\quad \prod_{p \le k} p \le 4^k.
\end{equation*}
For $n = 1$, this property holds because
\begin{equation*}
    \prod_{p \le 1} p = 1 \le 4^1 = 4,\qquad\prod_{p \le 2} p = 2 \le 4^2 = 16.
\end{equation*}
Now let $n \ge 1$, assume $\mathscr{P}_n$ is true, and prove $\mathscr{P}_{n + 1}$. Let $k \in \{1, \dots, 2n + 2\}$. We must show that $\prod_{p \le k} p \le 4^k$. If $k \le 2n$, this follows directly from $\mathscr{P}_n$. It remains to verify the inequality for $k = 2n + 1$ and $k = 2n + 2$. Since $2n + 2$ is not prime, we have
\begin{equation*}
    \prod_{p \le 2n + 2} p = \prod_{p \le 2n + 1} p,
\end{equation*}
so it suffices to prove it for $k = 2n + 1$. We can write
\begin{equation*}
    \prod_{p \le k} p = \prod_{p \le 2n + 1} p = \Bigl(\prod_{p \le n + 1} p\Bigr)\Bigl(\prod_{n + 1 < p \le 2n + 1} p\Bigr).
\end{equation*}
By property $\mathscr{P}_n$, the first factor is $\le 4^{n + 1}$, and we have seen that the second factor is $\le 4^n$.  Thus,
\begin{equation*}
    \prod_{p \le 2n + 1} p \le 4^{2n + 1} = 4^k.
\end{equation*}
This completes the induction. 

We have also the power series expansion
\begin{equation*}
    e^x = \sum_{k = 0}^{\infty} \frac{x^k}{k!}, \quad x \in \mathbb{R}.
\end{equation*}
It follows that for all $x > 0$ and all integers $k \ge 1$,
\begin{equation*}
    e^x > \frac{x^k}{k!},
\end{equation*}
which gives, in particular, for an integer $m \ge 1$,
\begin{equation*}
    e^m > \frac{m^m}{m!}.
\end{equation*}
Multiplying the latter inequality by $m!/e^m$ leads to
\begin{equation}\label{res1}
    m! > \Bigl(\frac{m}{e}\Bigr)^m.
\end{equation}

Let $n \ge 2$ be an integer. For all $i \ge 1$, let $p_i$ denote the $i$-th prime number. We have $i \le p_i$ for all $i$. Therefore,
\begin{equation*}
    \pi(n)! = 1 \cdot 2 \cdots \pi(n) \le p_1 p_2 \cdots p_{\pi(n)}.
\end{equation*}
By definition of $\pi(n)$, these primes are precisely those $\le n$. Hence
\begin{equation*}
    \pi(n)! \le \prod_{p \le n} p \le 4^n.
\end{equation*}
From Eq. (\ref{res1}), the left-hand side satisfies
\begin{equation*}
    \pi(n)! > \Bigl(\frac{\pi(n)}{e}\Bigr)^{\pi(n)},
\end{equation*}
so we obtain
\begin{equation*}
    \Bigl(\frac{\pi(n)}{e}\Bigr)^{\pi(n)} \le 4^n.
\end{equation*}
Taking natural logarithm of both sides gives
\begin{equation}\label{res2}
    \pi(n)(\log \pi(n) - 1) \le n\log 4.
\end{equation}
On $[1, +\infty)$, the function $x \mapsto x\log x - x$ is continuously differentiable and strictly increasing, since its derivative on $(1, +\infty)$ is
\begin{equation*}
    (\log x + 1) - 1 = \log x > 0.
\end{equation*}
Hence $x\log x - x$ is strictly increasing on $[1, +\infty)$. Suppose that there exists an integer $N \ge 3$ such that
\begin{equation*}
    \pi(N) > \frac{e N}{\log N}.
\end{equation*}
As the function $x \mapsto x\log x - x$ is increasing and $\frac{e N}{\log N} > 1$, we have
\begin{equation*}
    \frac{e N}{\log N}\Bigl(\log\Bigl(\frac{e N}{\log N}\Bigr) - 1\Bigr) < \pi(N)(\log \pi(N) - 1).
\end{equation*}
Using inequality (\ref{res2}), this implies
\begin{equation*}
    \frac{e N}{\log N}\Bigl(\log\Bigl(\frac{e N}{\log N}\Bigr) - 1\Bigr) < N\log 4.
\end{equation*}
Simplifying gives
\begin{equation*}
    e(1 + \log N - \log \log N - 1) < \log 4 \log N,
\end{equation*}
that is,
\begin{equation*}
    e\log N < \log 4 \log N + e\log \log N,
\end{equation*}
and hence
\begin{equation*}
    e < \log 4 + e\frac{\log \log N}{\log N},
\end{equation*}
or equivalently,
\begin{equation*}
    \frac{e-\log 4}{e} < \frac{\log \log N}{\log N}.
\end{equation*}

Let us now consider the function $x\longmapsto u(x)=\log(x)/x$ defined on the interval $[1, +\infty[$. The function $u$ is $\mathcal{C}^1$ on the interval and
\begin{equation*}
    u'(x)=\frac{1-\log(x)}{x}
\end{equation*}
which cancels for $x = e$. The second derivative of $u$ is given by
\begin{equation*}
    u''(x)=\frac{2\log x-3}{x^3},
\end{equation*}
and we get
\begin{equation*}
    u''(e)=-\frac{1}{e^3}<0,
\end{equation*}
meaning that the function $u$ has a global maximum for $x = e$. We have $u(e)=\log(e)/e = 1/e$ for all $x \ge 1$. We have
\begin{equation*}
    \frac{e-\log 4}{e} < \frac{\log \log N}{\log N} \le \frac{1}{e},
\end{equation*}

and therefore,
\begin{equation*}
    e-\log 4 < 1 \quad\Rightarrow\quad e < 1 + \log 4.
\end{equation*}
But $e \approx 2.71828$ and $1 + \log 4 \approx 2.38629$, which is a contradiction. Thus,
\begin{equation*}
    \pi(n) \le \frac{e n}{\log n}
\end{equation*}
for all integers $n \ge 3$.

\end{document}